\documentclass[11pt]{article}
\usepackage{amsmath,amssymb}

\def\r{\rightarrow}

\def\qed{\hfill\vrule height5pt width5pt depth0pt}
\def\one #1{1_{\{#1\}}}

\def\1{\mbox{\bf 1}}

\newcommand{\proof}{\noindent {\bf Proof:\ }}

\newtheorem{Theorem}{Theorem}
\newtheorem{Lemma}{Lemma}

\newtheorem{Corollary}{Corollary}
\newtheorem{Remark}{Remark}

\begin{document}
\title{Useful martingales for stochastic storage processes with L\'evy-type input}
% with application to decomposition results for L\'evy queues with service interruptions
\author{
Offer Kella\thanks{Department of Statistics; The Hebrew University of Jerusalem; Mount Scopus,
Jerusalem 91905; Israel ({\tt offer.kella@huji.ac.il})}
\thanks{Supported in part by grant 434/09 from the Israel Science
Foundation, the Vigevani Chair in Statistics and visitor grant No. 040.11.257 from The Netherlands Organisation for Scientific Research.}
\and
Onno Boxma\thanks{EURANDOM and Department of Mathematics and
Computer Science; Eindhoven University
of Technology; P.O. Box 513; 5600 MB Eindhoven; The Netherlands ({\tt boxma@win.tue.nl})}
}
\date{September 19, 2012}
\maketitle

\abstract{In this paper we generalize the martingale of Kella and Whitt to the setting of L\'evy-type processes and show that the (local) martingales obtained are in fact square integrable martingales which upon dividing by the time index converge to zero a.s. and in $L^2$. The reflected L\'evy-type process is considered as an example.
}

\vspace{0.1in}

\vspace{0.1in}


\section{Introduction}

In \cite{kw1992} a certain (local) martingale associated with L\'evy processes and its various applications is discussed (see also Section~IX.3 of \cite{a2003} and Section~4.4 of \cite{kip2006}). This has become a standard tool for studying various storage systems with L\'evy inputs and other problems associated with L\'evy process modeling. In \cite{ak2000} a generalization to a multidimensional (local) martingale associated with Markov additive processes with finite state space Markov modulation is considered, and in \cite{ak2001} a special case of the martingale of \cite{kw1992} for a reflected and a nonreflected L\'evy process with no negative jumps and applications to certain hitting times associated with these processes. A generalization to martingales associated with more general functions (than exponential) is given in \cite{ny2005}. The focus is on reflected and nonreflected processes but the main results seem to hold for the more general structure considered in \cite{kw1992}. There are many papers which apply this and related martingales. As these particular applications are not the scope of this study, we will not attempt to list them here.

The first goal of this paper is to extend the local martingale results of \cite{kw1992} to the case where the driving process is a L\'evy-type process. That is, it is a sum of stochastic integrals of some bounded left continuous right limit process with respect to coordinate processes associated with some multidimensional L\'evy process. Such processes with an even more general (predictable) integrand are discussed in \cite{applebaum}. The second goal is to extend the original results to show that without any further conditions the resulting local martingales are in fact square integrable martingales which upon division by the time parameter $t$, converge to zero almost surely and in $L^2$ as $t\r\infty$. Therefore, certain conditions originally made in \cite{kw1992} to ensure that the local martingales established there are martingales turn out to be unnecessary as Theorem~\ref{Mr0} here in particular applies to the special case ({\em L\'evy}, rather than {\em L\'evy-type}) treated there.

This article is organized as follows. In Section~\ref{sec:martingale} we develop the main local martingale.
In Section~\ref{sec:Mtt} we show that it is in fact a square integrable martingale and that its rate (defined appropriately) is zero almost surely and in $L^2$. In Section~\ref{sec:reflection} we give a small demonstration of the results with a reflected L\'evy-type process and a strong law for L\'evy-type processes which is established under some suitable assumption.

Although following the derivations requires some knowledge, we believe that the final results (in particular Theorem~\ref{Mr0} together with Theorem~\ref{martingale} and to a large extent also Corollary~\ref{cor:sl} and Theorem~\ref{thm:reflected}) are quite easy to use also by those who are not L\'evy process experts nor familiar with the theory of stochastic integration. One particular application that motivated this study is to establish decomposition results for L\'evy-driven polling systems (e.g., \cite{BIKM} and references therein), or more generally, on/off storage systems with L\'evy inputs, where the on/off structure could be quite general: during off times the process behaves like a subordinator and during on times it behaves like a (possibly unrelated) reflected L\'evy process. These results will be discussed in a separate paper which generalizes \cite{k1998,kw1991}, where the results established here are essentially needed and simplify the analysis considerably.

\section{A more general local martingale}
\label{sec:martingale}
For what follows given a c\`adl\`ag (right continuous left limit) function $g:\mathbb{R}_+\r\mathbb{R}$ we denote $g(t-)=\displaystyle \lim_{s\uparrow t}g(s)$, $\Delta g(t)=g(t)-g(t-)$ with the convention that $\Delta g(0)=g(0)$ and if $g$ is VF (finite variation on finite intervals), then $g^d(t)=\sum_{0\le s\le t}\Delta g(s)$ and $g^c(t)=g(t)-g^d(t)$. Also, $\mathbb{R}_+=[0,\infty)$, $\mathbb{R} =(-\infty,\infty)$  and {\em a.s.} abbreviates {\em almost surely}.

Let $X=(X_1,\ldots,X_K)$ be a c\`adl\`ag $K$-dimensional L\'evy process with respect to some standard filtration $\left\{\mathcal{F}_t|\ t\ge0\right\}$ with exponent
\begin{equation}
\psi(\alpha)=ic^T\alpha-\frac{\alpha^T\Sigma\alpha}{2}+\int_{\mathbb{R}^K}\left(e^{i\alpha^Tx}-1-i\alpha^Tx\one{\|x\|\le 1}\right)\nu(dx)
\end{equation}
where $^T$ denotes transposition, $\Sigma$ is positive semidefinite and $\|x\|=\sqrt{x^Tx}$.
When $X_1,\ldots,X_K$ have no negative jumps, then for any vector $\alpha\ge 0$ the Laplace-Stieltjes exponent is
\begin{eqnarray}\label{LSE}
\varphi(\alpha)&=&\log Ee^{-\alpha^TX(1)}=\psi(i\alpha)\\ &=&-c^T\alpha+\frac{\alpha^T\Sigma\alpha}{2}+\int_{\mathbb{R}_+^K}\left(e^{-\alpha^Tx}-1+\alpha^Tx\one{\|x\|\le 1}\right)\nu(dx)\ .\nonumber
\end{eqnarray}
It is well known that in this case $\varphi(\alpha)$ is finite for each $\alpha\ge 0$, that it is convex (thus continuous) with $\varphi(0)=0$, is infinitely differentiable in the interior of $\mathbb{R}_+$ and that for every $\alpha\ge 0$ for which $\alpha^TX$ is not a subordinator (not nondecreasing), $\varphi(t\alpha)\r\infty$ as $t\r\infty$. Furthermore, $EX_k(t)=-t\frac{\partial\varphi}{\partial\alpha_k}(0+)$ (finite or $+\infty$, but can never be $-\infty$) and when the first two right derivatives at zero are finite, then $\mbox{Cov}(X_k(t),X_\ell(t))=t\frac{\partial^2\varphi}{\partial\alpha_k\partial\alpha_\ell}\varphi(0+)$.

\begin{Lemma}
Let $I=(I_1,\ldots,I_K)$ be a bounded $K$-dimensional adapted c\`adl\`ag process. Then
\begin{equation}
e^{i\sum_{k=1}^K\int_{(0,t]}I_k(s-)dX_k(s)-\int_0^t\psi(I(s))ds}
\end{equation}
is a (complex valued) martingale.
When in addition $X_k$ have no negative jumps and $I_k$ are nonnegative then
\begin{equation}
e^{-\sum_{k=1}^K\int_{(0,t]}I_k(s-)dX_k(s)-\int_0^t\varphi(I(s))ds}
\end{equation}
is a real valued martingale.
\end{Lemma}

\proof
 Follows, for example, by applying a multidimensional generalization  of Corollary 5.2.2 and  Theorem 5.2.4 on pages 253-254 of \cite{applebaum} to the process
\begin{eqnarray}
dY(t)&=&\left(\sum_{k=1}^Kc_kI_k(s)-\varphi(I(s))\right)dt\nonumber \\ \nonumber \\&&+\sum_{k=1}^K\Big(I_k(t)dB_k(t)+I_k(t-)x\tilde{N}_k(dt,dx)\\ \nonumber \\&&\qquad\qquad\qquad\qquad\ +I_k(t-)xN_k(dt,dx)\Big)\nonumber
\end{eqnarray}
where $Y$, $B_k$, $N_k$ and $\tilde{N}_k$ are the notations from \cite{applebaum} with the obvious additional index $k$. Since we will not use these notations in this paper we only mention them briefly here. Moreover, $Y$ will soon be used for something else, in line with \cite{kw1992} and \cite{ak2000}.\qed

Setting $Z(t)=\sum_{k=1}^K \int_{(0,t]}I_k(s-)dX_k(s)+Y(t)$, the exact same proof from \cite{kw1992} can be employed to prove the following, where $a\wedge b=\min(a,b)$. We recall here that in \cite{kw1992} the driving process was some one-dimensional L\'evy process $X$ rather than $\sum_{k=1}^K \int_{(0,t]}I_k(s-)dX_k(s)$.

\begin{Theorem}\label{martingale}
Let $X=(X_1,\ldots,X_K)$ be a L\'evy process with exponent $\psi$ and,  when it has no negative jumps, Laplace-Stieltjes exponent $\varphi$. Let $I=(I_1,\ldots,I_K)$ be bounded c\`adl\`ag and adapted. Assume that $Y$ is c\`adl\`ag, VF (a.s.) and adapted. Then
\begin{eqnarray}
M(t)&=&\int_0^t\psi(I(s))e^{iZ(s)}ds+e^{iZ(0)}-e^{iZ(t)}+i\int_0^te^{iZ(s)}dY^c(s)\nonumber\\ \\
&&+\sum_{0<s\le t}e^{iZ(s)}\left(1-e^{-i\Delta Y(s)}\right)\nonumber
\end{eqnarray}
is a local martingale.
%and if the expected variation of $Y^c$ is finite and \[E\sum_{0<s\le t}|\Delta Y(s)|\wedge 1<\infty\ ,\] then it is a zero mean martingale.

\noindent
When $Z$ is bounded below, $X_k$ have no negative jumps and $I_k$ are nonnegative, then
\begin{eqnarray}\label{eq:LSTmart}
M(t)&=&\int_0^t\varphi( I(s))e^{- Z(s)}ds+e^{- Z(0)}-e^{- Z(t)}-\int_0^te^{- Z(s)}dY^c(s)\nonumber\\ \\
&&+\sum_{0<s\le t}e^{- Z(s)}\left(1-e^{ \Delta Y(s)}\right)\nonumber
\end{eqnarray}
is a local martingale.
 %and when the expected variation of $Y^c$ is finite, \[E\sum_{0<s\le t}|\Delta Y(s)|\wedge1<\infty\] and $Z$ is bounded below (in particular non-negative), then it is a zero mean martingale.
\end{Theorem}

We note that in \cite{kw1992} it was assumed that the expected number of jumps of $Y$ on finite intervals is finite in order for the local martingale to be a martingale. It is easy to show with the same proof that the weaker condition \[E\sum_{0<s\le t}|\Delta Y(s)|\wedge 1<\infty\ ,\] is sufficient. For example, if $Y$ is a subordinator (a nondecreasing L\'evy process) then it satisfies this condition. Nevertheless, as we will later show that these local martingales are in fact square integrable martingales, even this condition seems unnecessary. We also remark that the condition that $Z$ is bounded below is not really necessary for (\ref{eq:LSTmart}) to be a local martingale, but we will need it later to show that it is a square integrable martingale with rate zero, which is the more important result that we are aiming at.

It may seem more general to consider the multidimensional process defined via $Z_\ell(t)=\sum_{k=1}^K\int_{(0,t]}I_{\ell k}(s-)dX_k(s)+Y_\ell$, but we immediately see that the one-dimensional process
\begin{equation}
\sum_{\ell=1}^L\alpha_\ell Z_\ell(t)=\sum_{k=1}^K\int_{(0,t]}\sum_{\ell=1}^L\alpha_\ell I_{\ell k}(s-)dX_k(s)+\sum_{\ell=1}^K\alpha_\ell Y_\ell(t)
\end{equation}
has the same structure, resulting in the following (local) martingales
\begin{eqnarray}
M(t)&=&\int_0^t\psi(\alpha^TI(s))e^{i\alpha^TZ(s)}ds+e^{i\alpha^TZ(0)}-e^{i\alpha^TZ(t)}\nonumber\\ \\
&&+i\sum_{\ell=1}^L\alpha_\ell\int_0^te^{i\alpha^TZ(s)}dY_\ell^c(s)+\sum_{0<s\le t}e^{i\alpha^TZ(s)}\left(1-e^{-i\alpha^T\Delta Y(s)}\right)\nonumber
\end{eqnarray}
and
\begin{eqnarray}
M(t)&=&\int_0^t\varphi(\alpha^TI(s))e^{-\alpha^TZ(s)}ds+e^{-\alpha^TZ(0)}-e^{-\alpha^TZ(t)}\nonumber\\ \\
&&-\sum_{\ell=1}^L\alpha_\ell\int_0^te^{-\alpha^TZ(s)}dY_\ell^c(s)+\sum_{0<s\le t}e^{-\alpha^TZ(s)}\left(1-e^{\alpha^T\Delta Y(s)}\right)\nonumber
\end{eqnarray}
where $I$ is an $L\times K$-matrix valued function.

We note that when $J$ is a (right continuous) continuous time Markov chain with states $1,\ldots,K$, then with $I_k(t)=\one{J(t)=k}$ one has that $\sum_{k=1}^K\int_{(0,t]}I_k(s-)dX_k(s)$ is a Markov additive process. Adding additional jumps at state change epochs can be modeled by the process $Y$, which is obviously VF. For the case where $Y$ is continuous, this kind of a process and associated martingales were considered in \cite{ak2000}. The one-dimensional martingales considered here are not the same as the multidimensional ones considered there. However, the sum of the components of the latter does agree with the former.

We conclude this section with the following observation. Assume that $J$ is a c\`adl\`ag adapted process taking values in some finite set $1,\ldots,K$ (not necessarily Markovian). Let $I_k(t)=\alpha_k\one{J(t)=k}$. Then
\begin{equation}\psi(I(t))=\sum_{k=1}^K\psi_k(\alpha_k)\one{J(t)=k}\ ,
\end{equation}
where $\psi_k(\alpha_k)=\psi(0,...0,\alpha_k,0,\ldots,0)$ with $\alpha_k$ in the $k$th coordinate, is defined in the previous remark (and similarly with $\varphi$ when there are no negative jumps). Thus, in this case
\begin{equation}
\int_0^t\psi( I(s))e^{iZ(s)}ds=\sum_{k=1}^K\psi_k(\alpha_k)\int_0^te^{iZ(s)}\one{J(s)=k}ds\ .
\end{equation}
If in addition we replace $Y$ by $\beta Y$ for some $\beta\ge 0$ and denote $\tilde X_k(t)=\int_{(0,t]}\one{J(s)=k}dX_k(s)$ then
\begin{equation}Z(t)=\alpha^T\tilde X(t)+\beta Y(t)\end{equation}
and the (local) martingale becomes
\begin{eqnarray}
M(t)&=&\sum_{k=1}^K\psi_k(\alpha_k)\int_0^te^{iZ(s)}\one{J(s)=k}ds+e^{i Z(0)}-e^{iZ(t)}\nonumber\\ \\
&&+i\beta \int_0^te^{iZ(s)}dY^c(s)+\sum_{0<s\le t}e^{iZ(s)}\left(1-e^{-i\beta \Delta Y(s)}\right)\nonumber
\end{eqnarray}
and similarly
\begin{eqnarray}
M(t)&=&\sum_{k=1}^K\varphi_k(\alpha_k)\int_0^te^{-Z(s)}\one{J(s)=k}ds+e^{- Z(0)}-e^{-Z(t)}\nonumber\\ \\
&&-\beta \int_0^te^{-Z(s)}dY^c(s)+\sum_{0<s\le t}e^{-Z(s)}\left(1-e^{\beta \Delta Y(s)}\right)\nonumber
\end{eqnarray}
when there are no negative jumps.

It seems that the joint structure of $X$ is not important here. This is partly true in the sense that the evolution of the L\'evy part of the process during times when $J$ is at a given state is that of a one-dimensional L\'evy process. However, both $J$ and $Y$ may also depend on the joint structure.

\section{$M$ is a square integrable martingale with $M(t)/t\r0$ a.s. and in $L^2$\label{sec:Mtt}}
In this section we will show that $M$ is a square integrable martingale with $M(t)/t\r0$ a.s. and in $L^2$ as $t\to\infty$. This is something that was overlooked in \cite{kw1992}. To keep the discussion shorter, we will restrict it to the case of (\ref{eq:LSTmart}) where $Z$ is nonnegative, $X_k$ have no negative jumps and $I_k$ are nonnegative. The proofs for the general (complex valued) case are basically identical (but see Remark~\ref{fg}). Assuming the seemingly more general condition that $Z$ is bounded below rather than nonnegative is of no consequence to the proofs. We begin with the following.
\begin{Lemma}
Let $X$ be a semimartingale and $f\in\mathcal{C}^2$ (twice continuously differentiable). Denote by $[\cdot,\cdot]$ the quadratic variation process associated with a semimartingale. Then $f(X)$ is also a semimartingale with the following quadratic variation:
\begin{align}
[f(X),f(X)](t)&=\int_0^t\left(f'(X(s))\right)^2d[X,X]^c(s)\nonumber\\ &\quad+\sum_{0\le s\le t}\left(\Delta f(X(s))\right)^2 .
\end{align}
\end{Lemma}

\proof
Although this should have been a standard result in a book (such as \cite{p2004}) we did not find a direct reference. For its proof we apply the extended It\^o's Lemma (Thm. 32 on p. 78 of \cite{p2004}) to conclude that
\begin{eqnarray}
f(X(t))&=&f(X(0))+\int_{(0,t]}f'(X(s-))dX(s)\nonumber \\ \\ \nonumber&&+\mbox{continuous VF part}+\mbox{discrete VF part.}
\end{eqnarray}
As in the displayed equation following the definition of $[X,X]^c$ on p. 70 of \cite{p2004} we have that
\begin{equation}
[f(X),f(X)](t)=[f(X),f(X)]^c(t)+\sum_{0\le s\le t}\left(\Delta f(X(s))\right)^2 .
\end{equation}
Finally we note that the only term that can contribute to the continuous part of the quadratic variation associated with $f(X)$ is the stochastic integral part. Thus with the notation
$f(X_-)\cdot X(t)=\int_{(0,t]}f(X(s-))dX(s)$ we now have via Thm.~29 on p. 75 of \cite{p2004} that
\begin{eqnarray}
[f(X),f(X)]^c&=&[f(X_-)\cdot X,f(X_-)\cdot X]^c=\left(\left(f'(X_-)\right)^2\cdot[X,X]\right)^c\nonumber\\ \\ \nonumber&&=(\left(f'(X_-)\right)^2\cdot[X,X]^c
\end{eqnarray}
and the proof is complete.\qed

\begin{Corollary}\label{exp}
Assume that $X$ is a semimartingale, $Y$ is VF, adapted and $Z=X+Y$. Then
\begin{align}\label{e-Z}
[e^{-Z},e^{-Z}](t)&=\int_0^t e^{-2Z(s)}d[X,X]^c(s)\nonumber\\ &\quad +\sum_{0\le s\le t}e^{-2Z(s-)}\left(1-e^{-\Delta Z(s)}\right)^2 .
\end{align}
\end{Corollary}

\proof
Follows from $[Z,Z]^c=[X,X]^c$ (as $Y$ is VF), substitution and some obvious manipulations.\qed

\begin{Remark}\label{fg}\rm
Given the above, it is now an easy exercise to show that in fact for $X$ a semimartingale and  $f,g\in\mathcal{C}^2$ we have that
\begin{align}
[f(X),g(X)](t)&=\int_0^t f'(X(s))g'(X(s))d[X,X]^c(s)\nonumber\\ &\quad+\sum_{0\le s\le t}\Delta f(X(s))\Delta g(X(s))
\end{align}
and to conclude from this that, under the assumptions of Corollary~\ref{exp},
\begin{equation}\label{eiZ}
[e^{iZ},e^{iZ}](t)=\int_0^t e^{i2Z(s)}d[X,X]^c(s)+\sum_{0\le s\le t}e^{i2Z(s-)}\left(1-e^{i\Delta X(s)}\right)^2
\end{equation}
by treating the real and imaginary parts separately. This is needed for the general case which, as mentioned, is omitted from the discussion here.
\end{Remark}

Recall $M$ from (\ref{eq:LSTmart}) in Theorem~\ref{martingale} and that we are assuming that $X$ has no negative jumps and $I$ is nonnegative.

\begin{Corollary}\label{qv}
\begin{align}
[M,M](t)&=\int_0^t e^{-2Z(s)}d[\tilde X,\tilde X]^c(s)\nonumber\\ &\quad+\sum_{0< s\le t}e^{-2Z(s-)}\left(1-e^{-\Delta \tilde X(s)}\right)^2
\end{align}
\end{Corollary}

\proof The only part of $M$ that can contribute to the quadratic variation is
\begin{equation}
\sum_{0<s\le t}e^{-Z(s)}(1-e^{\Delta Y(s)})+e^{-Z(0)}-e^{-Z(t)}
\end{equation}
as the rest are continuous and VF. Clearly, only $e^{-Z(t)}$ contributes to the continuous part of this quadratic variation and, from Corollary~\ref{exp}, is given by $\int_0^te^{-2Z(s)}d[\tilde X,\tilde X]^c(s)$. Since $e^{-Z(0)}-e^{-Z(s)}=0$ for $s=0$ the `jump' at zero is excluded. Now, as $Z(s)=Z(s-)+\Delta \tilde X(s)+\Delta Y(s)$,
\begin{align}
&\Delta\left(\sum_{0<s\le t}e^{-Z(s)}(1-e^{\Delta Y})+e^{-Z(0)}-e^{-Z(t)}\right)(t)\nonumber\\
&=e^{-Z(t)}\left(1-e^{\Delta Y(t)}\right)+e^{-Z(t-)}-e^{-Z(t)}\nonumber\\
&=e^{-Z(t-)}\left(e^{-\Delta Y(t)}-1\right)e^{-\Delta\tilde X(t)}\\ &\qquad +e^{-Z(t-)}\left(1-e^{-\Delta Y(t)}e^{-\Delta \tilde X(t)}\right)\nonumber\\
&=-e^{-Z(t-)}\left(1-e^{-\Delta \tilde X(t)}\right)\ .\nonumber
\end{align}
As the discrete part of the quadratic variation is just the sum of squares of these jumps, we are done.\qed

\begin{Lemma}\label{MM}
\begin{eqnarray}
[M,M](t)&=&\int_0^t e^{-2Z(s)}A(s)ds+\tilde M(t)
\end{eqnarray}
where
\begin{equation}
A(s)=\varphi(2I(s))-2\varphi(I(s)) ,
\end{equation}
is nonnegative and bounded and $\tilde M$ is a martingale having bounded jumps.
\end{Lemma}

\proof
Recalling $\tilde X(t)=\sum_{k=1}^K\int_{(0,t]}I_k(s-)dX_k(s)$ we have
from Thm. 29 on p. 75 of \cite{p2004} that
\begin{eqnarray}
[\tilde X,\tilde X]=\sum_{k=1}^K\sum_{\ell=1}^K\left[I_k\cdot X_k,I_\ell\cdot X_\ell\right]=\sum_{k=1}^K\sum_{\ell=1}^KI_k I_\ell\cdot [X_k,X_\ell]
\end{eqnarray}
and thus also that
\begin{eqnarray}
[\tilde X,\tilde X]^c=\sum_{k=1}^K\sum_{\ell=1}^KI_k I_\ell\cdot [X_k,X_\ell]^c\ .
\end{eqnarray}
Now since we can write $X=B+C$, where $B$ is a Brownian motion and $C$ is a quadratic pure jump L\'evy process (e.g. see top of p. 71 of \cite{p2004}), then $[X_k,X_\ell]^c(t)=[B_k,B_\ell](t)=\sigma_{k\ell}t$ which implies that
\begin{align}\label{cont}
[\tilde X,\tilde X]^c(t)&=\int_0^t I(s)^T\Sigma I(s)ds\\ &=\int_0^t \left[\frac{(2I(s)^T)\Sigma (2I(s))}{2}-2 \frac{I(s)^T\Sigma I(s)}{2}\right]ds .
\nonumber\end{align}
Next, from $\Delta\tilde X(s)=\sum_{k=1}^KI(s-)\Delta X_k(s)$, we observe that
\begin{align}
&\sum_{0<s\le t}e^{-2Z(s-)}\left(1-e^{-\Delta\tilde X(s)}\right)\\&=\int_{(0,t]\times(0,\infty)^K}e^{-2Z(s-)}\left(1-e^{-I^T(s-)x}\right)N(ds,dx),
\nonumber\end{align}
where $N$ is the usual Poisson random measure with intensity measure $ds\otimes\nu(dx)$ associated with the jumps of $X$. Therefore, with $\tilde N(ds,dx)=N(ds,dx)-ds\otimes\nu(dx)$, recalling that $\int_{\mathbb{R}^K} (\|x\|^2\wedge 1) \nu(dx)<\infty$ and noting that $e^{-Z(s-)}(1-e^{-I(s-)x})\le (B\|x\|)\wedge 1$, where $B$ is an upper bound for $\|I(t)\|$ and thus $\int_0^t\int_{\mathbb{R}^K}Ee^{-2Z(s-)}(1-e^{-I(s-)x})^2\nu(dx)ds<\infty$, we have (e.g., Proposition~4.10 in \cite{r2004}) that
\begin{align}
\tilde M(t)=\int_{(0,t]\times(0,\infty)^K}e^{-2Z(s-)}\left(1-e^{-I^T(s-)x}\right)^2\tilde N(ds,dx),
\nonumber\end{align}
is a martingale, necessarily having bounded jumps, and so
\begin{align}\label{pmart}
&\sum_{0<s\le t}e^{-2Z(s-)}\left(1-e^{-\Delta\tilde X(s)}\right)^2\\&=\int_0^t\int_{(0,\infty)^K}e^{-2Z(s-)}\left(1-e^{-I^T(s-)x}\right)^2\nu(dx)ds+\tilde M(t)\ .
\nonumber\end{align}
%Next, we observe that since $X$ and thus its pure quadratic jump part is a L\'evy process then
%\begin{equation}
%\tilde N(t)=\sum_{0\le s\le t}\left(1-e^{-I(s-)^T\Delta X(s)}\right)^2-\int_0^t \left[\int_{(0,\infty)}\left(1-e^{-I(s)^Tx}\right)^2\nu(dx)\right]ds
%\end{equation}
%is a martingale having bounded jumps. To show this, one, e.g., first shows it for (multi-dimensional) compound Poisson processes and then takes appropriate limits. Thus also
%\begin{equation}
%\tilde M(t)=\int_{[0,t]}e^{-2Z(s-)}d\tilde N(s)\
%\end{equation}
%is a martingale (see Th. 51 on p. 38 and Th. 29 on p. 128 of \cite{p2004}).
Finally we observe that for any $a,x\in\mathbb{R}_+^K$
\begin{eqnarray}
\left(1-e^{-a^Tx}\right)^2&&=\left(e^{-(2a)^Tx}-1+(2a)^Tx\one{\|x\|\le 1}\right)\nonumber \\ \\ \nonumber &&-2\left(e^{-a^Tx}-1+a^Tx\one{\|x\|\le 1}\right)
\end{eqnarray}
and upon replacing $a$ by $I(s-)$ and integrating with respect to $\nu(dx)$, then together with (\ref{cont}) and (\ref{pmart}), the result is obtained.\qed

\begin{Theorem}\label{Mr0}
$M$ is a square integrable martingale with $M(t)/t\r 0$  as $t\r\infty$ a.s. and in $L^2$.
\end{Theorem}

\proof
 Since $I$ is bounded and $\varphi$ is continuous, then so is $\varphi(I)$. Therefore there exists a constant $C$ such that $\varphi(2I(s))-2\varphi(I(s))\le C$ and thus also $e^{-2Z(s)}(\varphi(2I(s))-2\varphi(I(s)))\le C$. As $\tilde M$ is a zero mean martingale, Lemma~\ref{MM} implies that
 \begin{equation}
 E[M,M](t)=\int_0^t e^{-2Z(s)}(\varphi(2I(s))-2\varphi(I(s)))ds\le Ct<\infty
 \end{equation}
  and thus, by Cor. 3 on p. 73 of \cite{p2004}, $M$ is a square integrable martingale with $EM^2(t)=E[M,M](t)$.
  Similarly, as $\int_{(0,t]}(1+s)^{-2}d\tilde M(s)$ is a zero mean martingale then
\begin{equation}\label{MMC}
E\int_0^t(1+s)^{-2}d[M,M](s)\le C\int_0^t(1+s)^{-2}ds=C\left(1-\frac{1}{1+t}\right)\le C .
\end{equation}
Letting $t\r\infty$ and applying monotone convergence on the left hand side (again with Cor. 3 on p. 73 of \cite{p2004}) implies that $\int_0^t(1+s)^{-1}dM(s)$ is a square integrable martingale with second moment given by the left side of (\ref{MMC}), that $\int_0^\infty(1+s)^{-1}dM(s)$ converges a.s. and thus, Ex. 14 on p. 95 of \cite{p2004} implies that $M(t)/(1+t)\r0$, hence also $M(t)/t\r0$ a.s.\qed

\section{A consequence for the reflected L\'evy-type process}\label{sec:reflection}
Reflected processes are widely used as models for various storage processes. With
\begin{equation}
L(t)=\displaystyle -\inf_{0\le s\le t}(Y(0)+\tilde X(s))^-
\end{equation}
it is well known that $Z(t)=\tilde X(t)+L(t)=0$ at any point of (right) increase of $L$ (e.g., \cite{k2006}). In the case where $X$ has no negative jumps $L$ is continuous. Thus in the general case $M$ becomes
\begin{eqnarray}
M(t)&=&\int_0^t\psi(I(s))e^{iZ(s)}ds+e^{iZ(0)}-e^{iZ(t)}+iL^c(t)\nonumber\\ \\
&&+\sum_{0<s\le t}\left(1-e^{-i\Delta L(s)}\right)\nonumber\\
&=&\int_0^t\psi(I(s))e^{iZ(s)}ds+e^{iZ(0)}-e^{iZ(t)}+iL(t)\nonumber\\ \\
&&-\sum_{0<s\le t}\left(e^{-i\Delta L(s)}-1+i\Delta L(s)\right)\nonumber
\end{eqnarray}
and when $X_k$ have no negative jumps and $I_k$ are nonnegative, then $\Delta L(s)=0$ and
\begin{eqnarray}
M(t)&=&\int_0^t\varphi( I(s))e^{- Z(s)}ds+e^{- Z(0)}-e^{- Z(t)}-L(t)\ .
\end{eqnarray}
By Theorem~\ref{Mr0} we therefore have for this case that
\begin{equation}\label{L(t)}
\frac{1}{t}\int_0^t\varphi(I(s))e^{-Z(s)}ds-\frac{1}{t}L(t)\r0
\end{equation}
a.s. and in $L^2$.

Also, we recall from the arguments of Theorem~1 of \cite{kw1996} that (path-wise) $\tilde X(t)/t\to \xi$ if and only if
\begin{equation}
\left(\frac{Z(t)}{t},\frac{L(t)}{t}\right)\to(\xi^+,-\xi^-)
\end{equation}
where $a^+=\max(a,0)$ and $a^-=\min(a,0)$. This is true for any c\`adl\`ag $\tilde X$, not necessarily having the special structure we consider here.

Thus, when $X_k$ have no negative jumps and $I_k$ are nonnegative it now follows that
\begin{equation}
\frac{1}{t}\int_0^t\varphi(I(s))e^{-Z(s)}ds\r -\xi^- .
\end{equation}
To figure out what $\xi$ is in this case, use the following result which is related to Theorem~\ref{Mr0} and holds regardless of whether there are negative jumps or not.
\begin{Lemma}\label{compensation}
Let $X$ be a one-dimensional L\'evy process with L\'evy measure $\nu$ satisfying \[\int_{|x|>1}|x|\nu(dx)<\infty\]  (equivalently $E|X(1)|<\infty$). Then for any bounded c\`adl\`ag adapted process $A$,
\begin{equation}\label{eq:compensation}
\frac{\int_{(0,t]}A(s-)dX(s)-EX(1)\int_0^t A(s)ds}{t}\r0
\end{equation}
a.s.
\end{Lemma}

\proof
Assume that $|A(t)|\le B<\infty$. Set, for $M>0$,
\begin{eqnarray}
X_M(t)&=&\sum_{0<s\le t}\Delta X(s)\one{\Delta X(s)>M}, \nonumber \\
X_{-M}(t)&=&\sum_{0<s\le t}\Delta X(s)\one{\Delta X(s)<-M}, \\
X_0(t)&=&X(t)-X_M(t)-X_{-M}(t). \nonumber
\end{eqnarray}
Also, denote $\xi_i=EX_i(1)$ for $i=M,-M,0$. Then $X_M, X_{-M}, X_0$ are independent L\'evy processes. $X_M$ is nondecreasing and $X_{-M}$ is nonincreasing. Now,
\begin{equation}
\left|\frac{1}{t}\int_{(0,t]}A(s-)dX_M(s)\right|\le B \frac{X_M(t)}{t}\ ,
\end{equation}
and by the strong law of large numbers for L\'evy processes we have that a.s.
\begin{equation}
\limsup_{t\r\infty}\left|\frac{1}{t}\int_{(0,t]}A(s-)dX_M(s)\right|\le B\xi_M=B\int_{(M,\infty)}x\nu(dx) .
\end{equation}
Clearly, we also have that
\begin{equation}
\left|\frac{1}{t}\int_0^t A(s)ds\right|\le B
\end{equation}
and thus
\begin{equation}
\limsup_{t\r\infty}\left|\frac{\int_{(0,t]}A(s-)dX_M(s)-\xi_M\int_0^tA(s)ds}{t}\right|\le 2B\int_{(M,\infty)}x\nu(dx)\ .
\end{equation}
Similarly
\begin{equation}
\limsup_{t\r\infty}\left|\frac{\int_{(0,t]}A(s-)dX_{-M}(s)-\xi_{-M}\int_0^tA(s)ds}{t}\right|\le 2B\int_{(-\infty,-M)}|x|\nu(dx)\ .
\end{equation}

Next, we observe that the martingale $M_0(t)=X_0(t)-\xi_0t$ is a L\'evy process with bounded jumps and thus its quadratic variation is a nondecreasing L\'evy process with bounded jumps which can also be compensated by a linear function to create a martingale (as in Lemma~\ref{MM}). Thus, like in the proof of Theorem~\ref{Mr0},  this implies that
\begin{equation}
\frac{\int_{(0,t]}A(s-)dX_0(s)-\xi_0\int_0^tA(s)ds}{t}\r 0
\end{equation}
a.s. (and also in $L^2$, but this is not needed here). To conclude,  denoting $\xi=\xi_0+\xi_M+\xi_{-M}=EX(1)$, we now clearly have that, a.s.,
\begin{equation}
\limsup_{t\r\infty}\left|\frac{\int_{(0,t]}A(s-)dX(s)-\xi\int_0^tA(s)ds}{t}\right|\le 2B\int_{(-\infty,-M)\cup(M,\infty)}|x|\nu(dx)
\end{equation}
and letting $M\r\infty$, recalling that $\int_{|x|>1}|x|\nu(dx)<\infty$, the proof is complete. \qed

\begin{Remark} \rm (relation with PASTA)
We note that if $E|X(1)|<\infty$ and $EX(1)\ne 0$ then since $X(t)/t\r EX(1)$ a.s., (\ref{eq:compensation}) is equivalent to
\begin{equation}
\frac{1}{X(t)}\int_{(0,t]}A(s-)dX(s)-\frac{1}{t}\int_0^t A(s)ds \r0 ,
\end{equation}
and thus $\frac{1}{X(t)}\int_{(0,t]}A(s-)dX(s)$ converges a.s. if and only if
$\frac{1}{t}\int_0^t A(s)ds$ does, and the limits coincide. When $X$ is a Poisson process, this is no less than an equivalent statement of the famous and often cited PASTA (Poisson Arrivals See Time Averages) property. See \cite{mw1990} for a martingale approach in a (nonexplosive) point process setting.
\end{Remark}
An immediate corollary of Lemma~\ref{compensation} is the following, where $\nu_k$ is the (marginal) L\'evy measure associated with $X_k$.
\begin{Corollary}\label{cor:sl} {\rm (strong law for $\tilde X$)}
Assume that  $\int_{|x|>1}|x|\nu_k(dx)<\infty$ (equivalently, $E|X_k(1)|<\infty$), for each $k$, and that
\begin{align}
\frac{1}{t}\int_0^tI_k(s)ds\r\beta_k
\end{align}
a.s., as $t\to\infty$. Then, a.s.,
\begin{align}
\xi=\lim_{t\to\infty}\frac{\tilde X(t)}{t}=\sum_{k=1}^K\beta_kEX_k(1) .
\end{align}
\end{Corollary}
Thus, we can summarize with the following.
\begin{Theorem}\label{thm:reflected}
Assume that, for each $k$, $X_k$ have no negative jumps, $\int_{x>1}x\nu_k(dx)<\infty$ (equivalently $EX_k(1)<\infty$) and that $I_k$ are nonnegative with
\begin{align}
\frac{1}{t}\int_0^tI_k(s)ds\r\beta_k
\end{align}
a.s., as $t\to\infty$. Then, a.s.,
\begin{align}
\frac{1}{t}\int_0^t\varphi(I(s))e^{-Z(s)}ds\to-\left(\sum_{k=1}^K\beta_kEX_k(1)\right)^-\ .
\end{align}
\end{Theorem}
Note that when $K=1$ and $I_1(t)=\alpha$ for all $t$, we have that $\beta_1=\alpha$ and that $\xi=EX_1(1)=-\varphi'(0)$. This immediately implies that if $\varphi'(0)>0$ then
\begin{align}\label{intez}
\frac{1}{t}\int_0^te^{-Z(s)}ds\to\frac{\alpha\varphi'(0)}{\varphi(\alpha)},\quad\mbox{as}\ t\to\infty\ .
\end{align}
When $\varphi'(0)<0$ or $\varphi'(0)=0$ but $X_1$ is not identically zero (so that $\varphi(\alpha)>0$ for each $\alpha>0$), then the limit is zero. When $\varphi'(0)>0$, this limit is the well known generalized Pollaczek-Khinchine formula. We also observe that $\varphi'(0)<0$ is the transient case and $\varphi'(0)=0$ but $X_1$ is not identically zero is the null recurrent case, so that neither is really a big surprise, but it is nice to see that it also follows directly from the above.

\end{document}